\documentclass[reqno]{amsart}
\usepackage{amsmath, amsthm, amssymb, amstext, latexsym, bbm}
\usepackage{hyperref,xcolor}
\hypersetup{pdfborder={0 0 0},colorlinks}
\usepackage{url}
\usepackage{euscript}
\usepackage[english]{babel}
\usepackage{graphicx}
\usepackage{color, xcolor}
\usepackage{enumitem}

\newtheorem{theorem}{Theorem}[section]
\newtheorem{lemma}{Lemma}[section]

\newtheorem{corollary}{Corollary}[section]
\theoremstyle{definition}
\newtheorem{definition}{Definition}[section]
\newtheorem{remark}{Remark}[section]
\newtheorem{example}{Example}[section]

\newcommand{\N}{\mathbb{N}}

\newcommand{\R}{\mathbb{R}}

\newcommand{\norm}[1]{\|#1\|}
\newcommand{\D}{\Delta}

\newcommand{\sumo}[1]{\sum_{k=1}^{#1}}
\renewcommand{\l}{\left}
\renewcommand{\r}{\right}

\numberwithin{equation}{section}

\title[Parametric indefinite discrete problem]{Two solutions for a parametric second-order discrete problem with indefinite coefficients} 

\author[E. Amoroso]{Eleonora Amoroso}
\address[E. Amoroso]{Department of Engineering, University of Messina, 98166 Messina, Italy.}
\email{eleonora.amoroso@unime.it}
\thanks{Corresponding author: Eleonora Amoroso, \texttt{eleonora.amoroso@unime.it}.}

\author[G. D'Aguì]{Giuseppina D'Aguì}
\address[G. D'Aguì]{Department of Engineering, University of Messina, 98166 Messina, Italy.}
\email{dagui@unime.it}


\begin{document}
	\begin{abstract}
		We investigate the existence of multiple solutions for a parametric	second-order discrete boundary value problem with indefinite coefficients. The nonlinear terms are independent and not necessarily generated by a single function. Using a variational approach and critical point theory, we prove the existence of at least two positive solutions for suitable values of the parameter. Finally, we provide examples that illustrate the applicability of the obtained results.
	\end{abstract}
	
	\maketitle
	
\section{Introduction} 

Difference equations have been widely investigated in recent years due to their strong connections with several applied fields such as numerical analysis, control theory, population dynamics and neural networks. For general overviews on difference equations and their applications, we refer to \cite{Agarwal-book, Elaydi, Kelley-Peterson-book}.

A large part of the literature is devoted to the study of existence and multiplicity of solutions for discrete problems by means of variational methods and critical point theory and we mention, among others, \cite{Agarwal-Perera-ORegan, Amo-Ca-Da, Amo-Ca-Ma, BonCandDAg, DAgui-Mawhin-Sciammetta} and the references therein. This approach allows one to translate the problem into the search of critical points of a suitable energy functional defined on a finite-dimensional Banach space.

In many contributions, the coefficients of the equation are assumed to be positive or nonnegative in order to guarantee the coercivity of the associated functional. However, from both theoretical and applied viewpoints, it is natural to consider problems with \textit{indefinite coefficients}, that is coefficients that may change sign. This situation arises, for instance, in models of neural network type, where inhibitory and excitatory effects coexist and lead to sign-changing potentials, see for instance \cite{Amoroso-discreto}.

The study of problems with indefinite coefficients has been recently addressed in the differential setting for the first time in \cite{Power}, where variational methods have been successfully applied. More recently, this approach has been extended to the discrete framework in \cite{Amoroso-discreto}, where a suitable equivalent norm was introduced in order to recover the variational structure of the problem even in the presence of possibly negative coefficients. In that paper, the existence of at least one solution was established without requiring global growth conditions on the nonlinear terms.

Motivated by these results, in the present paper we consider the parametric discrete problem described below and aim to obtain multiplicity results. In contrast with \cite{Amoroso-discreto}, where only the existence of one solution is proved, here we establish the existence of at least two positive solutions for the following problem
\begin{equation}\label{problem} \tag{P$_\textnormal{N}$}
	\begin{cases}
		-\Delta^2 u(k-1) + p(k)\,\Delta u(k-1) + q(k)\,u(k)
		= \lambda f_k(u(k)), & k\in[1,N],\\
		u(0) = u(N+1) = 0,
	\end{cases}
\end{equation}
where $N\in\mathbb{N}$, $[1,N]$ represents the discrete interval $\{1,2,\dots,N\}$ and $\lambda>0$ is a real parameter. Here, $\Delta u(k-1):=u(k)-u(k-1)$ denotes the forward difference operator, and consequently $\Delta^2 u(k-1)=u(k+1)-2u(k)+u(k-1)$ is the corresponding second-order difference operator. The functions $p,q:[1,N]\to\mathbb{R}$ are such that
\[
p(k)=\frac{s(k)}{s(k+1)}-1 \le 0, \qquad
\min_{k\in[1,N]} q(k) > -4\sin^2\!\left(\frac{\pi}{2(N+1)}\right)\frac{s(1)}{s(N+1)},
\]
where  $s:[0,N+2]\to\mathbb{R}$ is a positive and increasing function, that is, $s(k)>0$ for each $k\in[0,N+2]$ and $s(k)<s(k+1)$ for every $k\in[0,N+1]$. Note that, in general, the potential $q$ may change sign.  Finally, for each $k\in[1,N]$, the nonlinear terms $f_k:\mathbb{R}\to\mathbb{R}$ are assumed to be continuous and, in general, independent, so that they are not necessarily generated by a single function.

Our aim is to investigate the existence of multiple positive solutions to problem \eqref{problem} by means of variational methods. More precisely, by exploiting a two critical points theorem due to Bonanno and D'Agu\`{\i} \cite{Bonanno-DAgui-2016}, we prove the existence of at least two nontrivial critical points of the associated energy functional, corresponding to two positive solutions to \eqref{problem} with opposite energy sign. In particular, we provide suitable estimates ensuring the Palais--Smale condition and the unboundedness from below of the energy functional, even in the presence of indefinite coefficients. The multiplicity results are achieved under appropriate growth conditions on the nonlinear terms in an interval, or both at zero and at infinity.

The paper is organized as follows. In Section~\ref{Sec:prel}, we introduce the functional setting and recall the main tools from critical point theory. Section~\ref{Sec:main} is devoted to the proof of the main results, where we establish the existence of two positive solutions. Finally, we present some consequences and examples illustrating the applicability of our results.

\section{Preliminaries} \label{Sec:prel}
In this section we recall some basic definitions, properties and results that will be used throughout the paper. In particular, we introduce the variational framework, some properties of the associated energy functional and the abstract theorem that is the key tool of our investigation.

For any discrete interval $I_N$ and for any function $h \colon I_N \to \R$ we define
\[
h_- = \min_{k \in I_N} h(k), \qquad h_+ = \max_{k \in I_N} h(k). \]
In particular, although the weight $s$ is defined on $[0,N+2]$, only its values $s(1),\dots,s(N+1)$ appear in problem \eqref{problem} and in the functionals introduced below. Hence, we set
\[
s_- = \min_{k \in [1,N+1]} s(k), \qquad s_+ = \max_{k \in [1,N+1]} s(k).
\]
Moreover, for all $k \in [1,N]$ we set
\begin{align}\label{def:lambda_k}
	\mu_k=4\sin^2\l(\frac{k\pi }{2(N+1)}\r).
\end{align}
Consider the $N$-dimensional Banach space
\[
X=\{u \colon [0,N+1] \to \R \; \text{ such that } \; u(0)=u(N+1)=0\},
\]
endowed with the following norm, introduced recently by Amoroso et al. \cite{Amoroso-discreto},
\begin{equation*}
	\norm{u} =\l(\sum_{k=1}^{N+1}s(k)(\D u(k-1))^2+\sum_{k=1}^N s(k+1) q(k)(u(k))^2\r)^{\frac{1}{2}},
\end{equation*}
where we assume that $s:[0,N+2]\to\mathbb{R}$ is positive and increasing and $q:[1,N]\to\mathbb{R}$ is such that $q_- > -\mu_1\,\frac{s_-}{s_+}$. We introduce other standard norms in $X$ that will be used in the sequel
\begin{align*}
	&\norm{u}_\infty = \max_{k \in [1,N]} |u(k)|,\\
	&\norm{u}_2 = \l(\sumo{N}(u(k))^2 \r)^\frac{1}{2}, \\
	&\norm{u}_{\D,2} = \l(\sumo{N+1}(\D u(k-1))^2\r)^\frac{1}{2},
\end{align*}
for all $u \in X$. We know that they are all equivalent norms in $X$ and the following relations hold true, see \cite[p.~918]{BonCandDAg} for \eqref{equiv_bcd}-\eqref{norm_inf_bcd}, while \eqref{equiv_norm} can be obtained arguing as in \cite[Proposition 2.3]{Amoroso-discreto}, taking into account the estimate $s(k+1)q(k)\ge s_+q_-$ when $q_-<0$:
\begin{align}
	&\sqrt{\mu_1}\, \norm{u}_2\leq \norm{u}_{\D,2}\leq \sqrt{\mu_N}\, \norm{u}_2, 			 \label{equiv_bcd} \\[1ex]
	&\norm{u}_\infty\leq a_2\, \norm{u}_{\D,2},			\label{norm_inf_bcd}\\[1ex]
	&m\, \norm{u}_{\D,2}\leq \norm{u}\leq M\, \norm{u}_{\D,2}, 			\label{equiv_norm} 
\end{align}
for all $u \in X$, with $\mu_1, \mu_N$ defined in \eqref{def:lambda_k},
\begin{equation} \label{def:a_2}
	a_2=\begin{cases}
		\frac{\sqrt{N+1}}{2} &\text{if $N$ is odd},\\
		\l(\frac{2}{N}+\frac{2}{N+1}\r)^{-\frac{1}{2}}\quad &\text{if $N$ is even},
	\end{cases}
\end{equation}
and
\begin{equation} \label{def:mM}
	\begin{aligned}
		m=\begin{cases}
			\sqrt{s_-} & \text{if } q_-\geq 0,\\
			\sqrt{s_-\l(1+\frac{s_+}{s_-}\,\frac{q_-}{\mu_1}\r)}\; & \text{if } q_-<0,
		\end{cases} \quad
		M=\begin{cases}
			\sqrt{s_+\l(1+\frac{q_+}{\mu_1}\r)}\; & \text{if } q_+\geq 0,\\
			\sqrt{s_+} & \text{if }  q_+<0.
		\end{cases}
	\end{aligned}
\end{equation}
Notice that $a_2 \le \frac{\sqrt{N+1}}{2}$ for any $N \in \N$ and, since $s$ is positive and increasing, $s_- = s(1) > 0$ and $s_+ = s(N+1) >0$. Moreover, the assumption $q_- > -\mu_1\,\frac{s_-}{s_+}$, which reduces to $q_- > -\mu_1$ when $s$ is constant, guarantees that $m>0$, so that, by \eqref{equiv_norm}, $\norm{\cdot}$ is a norm. 

Our approach in studying the existence of solutions to problem \eqref{problem} is variational and it is based on critical point theory. In order to present the variational framework of our investigation, we assume that $f_k:\mathbb{R}\to\mathbb{R}$ is continuous for all $k \in [1,N]$. Put
\begin{equation*}
	F_k(t)=\int_0^t{f_k(\xi)\, d\xi}, \qquad \forall t\in \R, \quad \forall k\in [1,N],
\end{equation*}
and consider the so-called energy functional $I_\lambda \colon X \to \R$ related to problem \eqref{problem} defined as
\[
I_\lambda = \Phi - \lambda \Psi \qquad \text{for any } \lambda>0,
\]
where the functionals $\Phi, \Psi : X  \to \R$ are given by
\begin{align*}
	\Phi(u)=\frac{1}{2}\norm{u}^2, \qquad 
	\Psi(u)=\sum_{k=1}^N s(k+1)F_k(u(k)),
\end{align*}
for all $u \in X$. The functionals $\Phi,\Psi, I_\lambda$ are of class $C^1(X,\R)$, with 
\begin{align*}
	\Phi'(u)(v)&= \sumo{N+1} s(k) \Delta u(k-1) \Delta v(k-1) + \sumo{N} s(k+1) q(k) u(k) v(k) , \\
	\Psi'(u)(v) & =\sum_{k=1}^N s(k+1) f_k(u(k)) \, v(k), \\
	I_\lambda'(u)(v) &=  \Phi'(u)(v) - \lambda \Psi'(u)(v),
\end{align*}
for every $u,v\in X$. In \cite[Section 2]{Amoroso-discreto} the authors prove that $u \in X$ is a solution to problem \eqref{problem} if and only if it is a critical point of the energy functional $I_\lambda$, that is $I'_\lambda(u)(v)=0$ for all $v\in X$. Hence, establishing existence results for problem \eqref{problem} is equivalent to providing results on the existence of critical points of the energy functional $I_\lambda$. To this aim, we exploit a two critical points theorem presented by Bonanno-D'Aguì \cite[Theorem 2.1]{Bonanno-DAgui-2016} that we recall in the following together with a compactness condition exploited in the statement of the abstract result.

\begin{definition}
	We say that $I_\lambda$ satisfies the Palais-Smale condition (in short \textnormal{(PS)}-condition), if any sequence $\{u_n\}\subseteq X$ such that
	\begin{enumerate}[label=\textnormal{(PS$_{\arabic*}$)},ref=\textnormal{PS$_{\arabic*}$}]
		\item\label{P_1}$I_\lambda(u_n)$ is bounded,
		\medskip
		\item\label{P_2} $I_\lambda'(u_n) \to 0$ in ${X^\ast}$ as $n \to \infty$,
	\end{enumerate}
	has a strongly  convergent subsequence in $X$.
\end{definition}

\begin{theorem} \label{th:ZAA}
	Let $X$ be a real Banach space and let $\Phi, \Psi : X \rightarrow \R$ be two continuously G\^ateaux differentiable functionals such that $\inf_X \Phi = \Phi(0) = \Psi(0)=0$.
	
	Assume that there are $r \in \R$ and $\tilde{u} \in X$, with $0< \Phi(\tilde{u}) < r$, such that
	\begin{equation} \label{hp:th-zaa}
		\dfrac{\sup\limits_{u \in \Phi^{-1} (]- \infty, r])} \Psi(u)}{r} < \dfrac{\Psi(\tilde{u})}{\Phi(\tilde{u})},
	\end{equation}
	$I_\lambda = \Phi-\lambda\Psi$ satisfies \textnormal{(PS)}-condition and it is unbounded from below for all $\lambda\in {\Lambda}_r$, with
	\begin{align*}
		\Lambda_r = \l( \frac{\Phi(\tilde{u})}{\Psi(\tilde{u})},\frac{r}{\sup\limits_{u\in \Phi^{-1}(]-\infty,r])}\Psi(u)} \r).
	\end{align*}
	
	Then, for each $\lambda\in \Lambda_r$, the functional $I_\lambda=\Phi-\lambda\Psi$ admits at least two non-zero critical points $u_{\lambda,1}, u_{\lambda,2} \in X$ such that $I_\lambda (u_{\lambda,1}) < 0 < I_\lambda (u_{\lambda,2})$.
\end{theorem}

\section{Main results} \label{Sec:main}
In this section we present our main existence results for problem \eqref{problem}. We first establish some properties of the associated energy functional and then apply a two critical points theorem to obtain the existence of two positive solutions.

Throughout the paper, we assume the following conditions:
\begin{enumerate}[label=\textnormal{\textbf{(H)}},ref=\textnormal{H}]
	\item\label{H} 
	\begin{enumerate}[itemsep=0.2cm, label=$(\textnormal{h}_\arabic*)$, ref = $\textnormal{h}_\arabic*$]
		\item\label{hp:s} $s:[0,N+2]\to\mathbb{R}$ is positive and increasing;
		\item\label{hp:q-} $q:[1,N]\to\mathbb{R}$ is such that $q_- > -\mu_1\,\dfrac{s_-}{s_+}$;
		\item\label{hp:f_cont} $f_k:\mathbb{R}\to\mathbb{R}$ is continuous for all $k \in [1,N]$;
		\item\label{hp:f_sign} $f_k(0) \ge 0$ for all $k \in [1,N]$.
	\end{enumerate}
\end{enumerate}

\begin{remark}
	Since we are interested in positive solutions, we can truncate the nonlinearities before zero without any loss of generality. So, from now on, we consider $f_k(t)=f_k(0)$  for all $t<0, k \in [1,N]$.
\end{remark}

The following result can be found in \cite[Lemma 2.9]{Amoroso-discreto} and it provides information on the sign of the solution. 

\begin{lemma} \label{lemma:positive_sol}
	Let \eqref{H} be satisfied.
	If $u$ is a non-zero solution of problem \eqref{problem}, then $u$ is positive in $[1,N]$.
\end{lemma}

We first prove a result on the properties of the energy functional which will be used to apply Theorem~\ref{th:ZAA}.

\begin{lemma} \label{lemma:funzionale}
	Let \eqref{H} be satisfied and assume that
	\begin{enumerate}[itemsep=0.2cm, label=$(\textnormal{h}_\arabic*)$, ref = $\textnormal{h}_\arabic*$]
		\setcounter{enumi}{4}
		\item\label{hp:Linf} $L_\infty := \displaystyle\min_{k\in[1,N]} \liminf_{t \to + \infty} \dfrac{F_k(t)}{t^2} >0$.
	\end{enumerate}
	Then $I_\lambda$ satisfies the \textnormal{(PS)}-condition and it is unbounded from below for any $\lambda \in \l( \dfrac{M^2 \mu_N}{2 L_\infty s_-} , + \infty\r)$.
\end{lemma}
\begin{proof}
	First, we prove that $I_\lambda$ satisfies the (PS)-condition. So, let $\{u_n\}_{n \in \N} \subseteq X$ be such that \eqref{P_1}-\eqref{P_2} hold. We aim to prove that there exists a strongly convergent subsequence in $X$. For all $n \in \N$ consider
	\[
	u_n^+= \max \l\{ u_n, 0 \r\} \quad \text{ and } u_n^-= \max \l\{ -u_n, 0 \r\}.
	\]
	\textit{Step 1: $\{u_n^-\}$ is bounded for all $\lambda>0$.}\\
	
	If $u_n^-=0$ for $n$ large enough there is nothing to prove; otherwise, from \eqref{P_2} it follows that
	\begin{equation} \label{dim:step1}
		\lim\limits_{n \to + \infty} \dfrac{I'_\lambda(u_n)(u_n^-)}{\norm{u_n^-}} =0.
	\end{equation}
	Fix $n \in \N$ and $\lambda>0$. It is easy to prove that
	\begin{align*}
		\D u_n(k-1) \D u_n^-(k-1) \le - \l( \D u_n^-(k-1) \r)^2.
	\end{align*}
	Then, we have
	\begin{align*}
		\Phi'(u_n)(u_n^-)&= \sumo{N+1} s(k) \Delta u_n(k-1) \Delta u_n^-(k-1) + \sumo{N} s(k+1) q(k) u_n(k) u_n^-(k)\\
		&\le - \sumo{N+1} s(k) \l( \D u_n^-(k-1) \r)^2 - \sumo{N} s(k+1) q(k) \l(u_n^-(k)\r)^2\\
		&= - \norm{u_n^-}^2.
	\end{align*}
	Moreover, combining the definition of $u_n^-$ together with \eqref{hp:f_sign}, we obtain
	\begin{align*}
		\Psi'(u_n)(u_n^-) = \sumo{N} s(k+1) f_k(u_n(k)) u_n^-(k) = \sumo{N} s(k+1) f_k(0) u_n^-(k) \ge 0.
	\end{align*}
	Therefore, we get
	\begin{align*}
		I'_\lambda(u_n)(u_n^-) = \Phi'(u_n)(u_n^-) - \lambda \Psi'(u_n)(u_n^-) \le - \norm{u_n^-}^2,
	\end{align*}
	for all $n \in \N$. Exploiting the latter in \eqref{dim:step1} one has
	\begin{align*}
		\lim\limits_{n \to + \infty} \norm{u_n^-}=0,
	\end{align*}
	so that $\{u_n^-\}$ is bounded for every $\lambda>0$, namely there exists $M_1>0$ such that $\norm{u_n^-} \le M_1 $. Furthermore, using \eqref{equiv_bcd} and \eqref{equiv_norm}, we have
	\begin{align*} 
		\norm{u_n^-}_2 \le \dfrac{\norm{u_n^-}_{\D,2}}{\sqrt{\mu_1}} \le \dfrac{M_1}{m \sqrt{\mu_1}} = M_2,
	\end{align*}
	which implies that
	\begin{equation}\label{dim:step1-concl}
		u_n(k) \ge - M_2 \qquad \text{for all } n \in \N.
	\end{equation}
	
	\textit{Step 2: $\{u_n^+\}$ is bounded for all $\lambda>\dfrac{M^2 \mu_N}{2 L_\infty s_-}$.}\\
	
	Fix $\lambda>\dfrac{M^2 \mu_N}{2 L_\infty s_-}$, then there is $\ell>0$ such that
	\begin{align*}
		L_\infty > \ell > \dfrac{M^2 \mu_N}{2 s_- \lambda}.
	\end{align*}
	Arguing by contradiction, suppose that $\{u_n^+\}$ is unbounded. From the assumption $L_\infty>0$ it follows that
	\begin{align*}
		\liminf_{t \to + \infty} \dfrac{F_k(t)}{t^2} \ge L_\infty > \ell > \dfrac{M^2 \mu_N}{2 s_- \lambda} \qquad \text{for all } k \in [1,N].
	\end{align*}
	Hence, for all $k \in [1,N]$ there exists $\delta_k >0$ such that
	\begin{align} \label{dim:step2-a}
		F_k(t) > \ell t^2 \qquad \text{for all } t > \delta_k.
	\end{align}
	On the other side, for all $k \in [1,N], t \in [-M_2 , \delta_k]$ we have
	\begin{align*}
		F_k(t) &\ge \min_{t \in [-M_2 , \delta_k]} F_k(t) + \ell t^2 - \ell t^2\\
		&\ge \min_{t \in [-M_2 , \delta_k]} F_k(t) + \ell t^2 - \ell \l( \max\l\{ M_2, \delta_k \r\} \r)^2\\
		&\ge \ell t^2 - \max\l\{  \ell \l( \max\l\{ M_2, \delta_k \r\} \r)^2 - \min_{t \in [-M_2 , \delta_k]} F_k(t) , 0 \r\},
	\end{align*}
	that can be read as
	\begin{align} \label{dim:step2-b}
		F_k(t) \ge \ell t^2 - M_3(k) \qquad \text{for all } t \in [-M_2 , \delta_k],
	\end{align}
	with $M_3(k)\ge0$. Putting together \eqref{dim:step2-a} and \eqref{dim:step2-b}, we obtain
	\begin{align*}
		F_k(t) \ge \ell t^2 -M_3(k) \qquad \text{for all } t \ge -M_2, \; k \in [1,N].
	\end{align*}
	By step 1 (see also \eqref{dim:step1-concl}), we get
	\begin{align*}
		F_k(u_n(k)) \ge \ell (u_n(k))^2 - M_3(k) \qquad \text{for all } n \in \N, \; k \in [1,N].
	\end{align*}
	Then, it follows that
	\begin{align*}
		\Psi(u_n) &\ge \sumo{N} s(k+1) \l( \ell (u_n(k))^2 - M_3(k) \r)\\
		&\ge \ell s_- \norm{u_n}_2^2 - s_+ \bar{M}_3,
	\end{align*}
	for all $n \in \N$, with $ \bar{M}_3 = \sumo{N} M_3(k)$. Therefore, exploiting the latter together with \eqref{equiv_bcd} and \eqref{equiv_norm}, we obtain
	\begin{align*}
		I_\lambda(u_n) \le \l( \dfrac{1}{2} M^2 \mu_N - \lambda \ell s_- \r) \norm{u_n}^2_2 + \lambda s_+ \bar{M}_3 \qquad \text{for all } n \in \N.
	\end{align*}
	Since $\dfrac{1}{2} M^2 \mu_N - \lambda \ell s_- <0$ and $\norm{u_n}_2 \to + \infty$, one has $\lim\limits_{n \to + \infty} I_\lambda(u_n) = - \infty$ which is in contradiction with \eqref{P_1}.\\
	
	Hence, from step 1 and step 2 it follows that $\{u_n\}$ is bounded. Since $X$ is finite dimensional, $\{u_n\}$ admits a convergent subsequence, so that $I_\lambda$ satisfies (PS)-condition for all $\lambda \in \l( \dfrac{M^2 \mu_N}{2 L_\infty s_-} , + \infty\r)$.\\
	
	Finally, in order to prove that $I_\lambda$ is unbounded from below for any $\lambda > \dfrac{M^2 \mu_N}{2 L_\infty s_-}$ it is sufficient to consider a sequence $\{u_n\} \subset X$ such that $\{u_n^-\}$ is bounded and $\{u_n^+\}$ is unbounded. Then, arguing as before, we obtain that $\lim\limits_{n \to + \infty} I_\lambda(u_n) = - \infty$ and this completes the proof.
\end{proof}

We now state our main result on the existence of two positive solutions for problem \eqref{problem}. Set
\begin{align*}
	\sigma = s(1)+s(N+1)+\sumo{N} s(k+1)\,q(k), \qquad K = \dfrac{m^2}{a_2^2}\, \dfrac{s_-}{s_+}\, \dfrac{1}{\sigma},
\end{align*}
where $a_2$ is defined in \eqref{def:a_2}, $m,M$ are given by \eqref{def:mM}, and $s_-,s_+$ are defined in Section~\ref{Sec:prel}.
Notice that $\sigma>0$. Indeed, if $e \in X$ is such that $e(k)=1$ for all $k \in [1,N]$, a direct computation shows that $\sigma=\norm{e}^2>0$.
\begin{theorem} \label{th:main1}
	Let \eqref{H} be satisfied and assume that
	\begin{enumerate}[itemsep=0.2cm, label=$(\textnormal{h}_\arabic*)$, ref = $\textnormal{h}_\arabic*$]
		\setcounter{enumi}{5}
		\item\label{hp:cond1} there exist two constants $c,d$, with $0<d<c$ and $F_k(d)\ge0$ for all $k\in[1,N]$, such that
		\begin{align*}
			\dfrac{\displaystyle\sumo{N} \max_{|t| \le c} F_k(t)}{c^2} < K \min\l\{ \dfrac{\displaystyle\sumo{N} F_k(d)}{d^2} , \dfrac{\sigma\, L_\infty}{M^2 \mu_N} \r\}.
		\end{align*}
	\end{enumerate}
	Then, for every $\lambda \in \Lambda$, with
	\begin{align*}
		\Lambda = \left( \max\l\{ \dfrac{M^2 \mu_N}{2 L_\infty s_-} , \dfrac{\sigma}{2s_-} \dfrac{d^2}{\displaystyle\sumo{N} F_k(d)} \r\} , \dfrac{m^2}{2a_2^2 s_+}  \dfrac{c^2}{\displaystyle\sumo{N} \max_{|t| \le c} F_k(t)} \right),
	\end{align*}
	with the convention $\frac{1}{L_\infty}=0$ if $L_\infty=+\infty$, problem \eqref{problem} admits at least two positive solutions $u_1, u_2 \in X$ with opposite energy sign, i.e.
	\[
	u_1, u_2 >0 \text{ in } [1,N] \qquad \text{and} \qquad I_\lambda(u_1)<0<I_\lambda(u_2).
	\]
\end{theorem}
\begin{proof}
	We want to apply Theorem \ref{th:ZAA}. Choose $r= \dfrac{1}{2} \l( \dfrac{c \, m}{a_2} \r)^2$. If $\Phi(u)\le r$, then $\norm{u}\le \frac{cm}{a_2}$ and, by \eqref{norm_inf_bcd} and \eqref{equiv_norm},
	\begin{align*}
		\norm{u}_\infty \le a_2\norm{u}_{\D,2}\le \frac{a_2}{m}\norm{u}\le c \qquad \text{ for all } u \in \Phi^{-1}(]-\infty,r]).
	\end{align*}
	Since $\max_{|t|\le c}F_k(t)\ge F_k(0)=0$ for all $k\in[1,N]$, it follows that
	\begin{align*}
		\dfrac{\sup\limits_{u \in \Phi^{-1}(]-\infty,r])} \Psi(u)}{r} \le \dfrac{2 a_2^2 s_+}{m^2} \dfrac{\displaystyle\sumo{N} \max_{|t| \le c} F_k(t)}{c^2}.
	\end{align*}
	Let $\tilde{u} \in \R^{N+2}$ be defined as follows
	\begin{align*}
		\tilde{u}(0)=\tilde{u}(N+1)=0 \qquad \text{ and } \qquad \tilde{u}(k)=d \quad \text{for all } k \in [1,N].
	\end{align*}
	Clearly, $\tilde{u} \in X$ and
	\begin{align*}
		\Phi(\tilde{u}) = \dfrac{1}{2} \norm{\tilde{u}}^2 = \dfrac{1}{2} \l( s(1) d^2 + s(N+1) d^2 + \sumo{N} s(k+1) q(k) d^2 \r) = \dfrac{\sigma}{2}\, d^2>0,
	\end{align*}
	and, since $F_k(d)\ge0$ for all $k\in[1,N]$,
	\begin{align*}
		\Psi(\tilde{u})=\sumo{N}s(k+1)F_k(d)\ge s_- \sumo{N} F_k(d).
	\end{align*}
	Then, exploiting assumption \eqref{hp:cond1} and the definition of $K$, we obtain
	\begin{align*}
		\dfrac{\sup\limits_{u \in \Phi^{-1}(]-\infty,r])} \Psi(u)}{r} &\le \dfrac{2 a_2^2 s_+}{m^2} \dfrac{\displaystyle\sumo{N} \max_{|t| \le c} F_k(t)}{c^2} \\
		& < \dfrac{2 s_-}{\sigma} \dfrac{\displaystyle\sumo{N} F_k(d)}{d^2} \le \dfrac{\Psi(\tilde{u})}{\Phi(\tilde{u})}.
	\end{align*}
	Moreover, $\Phi(\tilde{u})<r$, that is $d^2<\frac{m^2}{a_2^2\sigma}\,c^2$. Indeed, arguing by contradiction, if $d^2 \ge \frac{m^2}{a_2^2\sigma}\,c^2$, then, since $d<c$ and $s_-\le s_+$,
	\begin{align*}
		\dfrac{\displaystyle\sumo{N} \max_{|t| \le c} F_k(t)}{c^2} \ge \dfrac{\displaystyle\sumo{N} F_k(d)}{c^2} \ge \dfrac{m^2}{a_2^2\sigma} \dfrac{\displaystyle\sumo{N} F_k(d)}{d^2} \ge K \dfrac{\displaystyle\sumo{N} F_k(d)}{d^2},
	\end{align*}
	against \eqref{hp:cond1}. Hence, hypothesis \eqref{hp:th-zaa} of Theorem \ref{th:ZAA} holds.
	From \eqref{hp:cond1}, since $\sum_{k=1}^N \max_{|t|\le c} F_k(t) \ge 0$ and $K>0$, we have
	\[
	\min\left\{
	\frac{\sum_{k=1}^N F_k(d)}{d^2},
	\frac{\sigma L_\infty}{M^2\mu_N}
	\right\} > 0.
	\]
	In particular, $L_\infty > 0$. Moreover, a direct computation shows that \eqref{hp:cond1} is equivalent to $\Lambda\neq\emptyset$, and $\Lambda\subseteq\Lambda_r$. Hence, Lemma \ref{lemma:funzionale} ensures that the energy functional $I_\lambda$ satisfies (PS)-condition and it is unbounded from below for all $\lambda \in \Lambda \subset \l( \dfrac{M^2 \mu_N}{2 L_\infty s_-} , + \infty\r)$. Then, by Theorem \ref{th:ZAA} it follows that  for all $\lambda \in \Lambda$ the energy functional $I_\lambda$ admits at least two critical points $u_1, u_2 \in X$ such that $I_\lambda(u_1)<0<I_\lambda(u_2)$, that is problem \eqref{problem} admits at least two solutions in $X$ with opposite energy sign. Finally, Lemma \ref{lemma:positive_sol} guarantees that the solutions are positive in $[1,N]$ and this completes the proof.
\end{proof}

\begin{remark}
	If the functions $f_k$ are nonnegative for all $k \in [1,N], t \ge 0$, then the functions $F_k$ are all nondecreasing and, owing to the truncation before zero, one has $\max_{|t| \le c} F_k(t)=F_k(c)$ for every $k \in [1,N]$.
\end{remark}

Now, we present several consequences of the main results under suitable growth conditions on the nonlinear terms both at zero and at infinity.

\begin{corollary} \label{corollary:lim1}
	Let \eqref{H} be satisfied and assume that
	\begin{enumerate}[itemsep=0.2cm, label=$(\textnormal{h}_\arabic*)$, ref = $\textnormal{h}_\arabic*$]
		\setcounter{enumi}{6}
		\item\label{hp:limitea0} $\limsup\limits_{t\to 0^+}
		\dfrac{\displaystyle\sum_{k=1}^N F_k(t)}{t^2}=+\infty$, and there exists $\bar t>0$ such that $F_k(t)\ge0$ for all $t\in[0,\bar t]$ and $k\in[1,N]$,
		\item\label{hp:limiteainfinito} $\liminf\limits_{t\to+\infty}
		\dfrac{F_k(t)}{t^2}=+\infty$ \; for all $k \in [1,N]$.
	\end{enumerate}
	Then, for every $\lambda \in (0, \lambda^*_1)$, with
	\begin{align*}
		\lambda^*_1	= \dfrac{m^2}{2 a_2^2 s_+} \, \sup_{c>0}	\dfrac{c^2}{\displaystyle\sumo{N} \max_{|t|\le c} F_k(t)},
	\end{align*}
	problem \eqref{problem} admits at least two positive solutions $u_1, u_2 \in X$ with opposite energy sign.
\end{corollary}
\begin{proof}
	By assumption \eqref{hp:limiteainfinito}, one has $L_\infty=+\infty$. Hence, the (PS)-condition and the unboundedness from below of the functional $I_\lambda$ follow from Lemma \ref{lemma:funzionale} for every $\lambda>0$. Now, fix $\lambda \in (0, \lambda^*_1)$ and $c>0$ such that
	\[
	\lambda <  \dfrac{m^2}{2 a_2^2 s_+} \, \dfrac{c^2}{\displaystyle\sumo{N} \max_{|t|\le c} F_k(t)}.
	\]
	From hypothesis \eqref{hp:limitea0}, there exists $d\in(0,\min\{c,\bar t\})$, so that $F_k(d)\ge0$ for all $k\in[1,N]$, such that
	\[
	\dfrac{2s_-}{\sigma} \dfrac{\displaystyle\sum_{k=1}^N F_k(d)}{d^2}
	> \dfrac{1}{\lambda} > \dfrac{2a_2^2 s_+}{m^2} \, \dfrac{\displaystyle\sumo{N} \max_{|t|\le c} F_k(t)}{c^2}.
	\]
	Moreover, since $L_\infty=+\infty$, we have
	\[
	\min\l\{ \dfrac{\displaystyle\sumo{N} F_k(d)}{d^2} , \dfrac{\sigma L_\infty}{M^2 \mu_N} \r\}
	= \dfrac{\displaystyle\sum_{k=1}^N F_k(d)}{d^2}.
	\]
	Therefore, condition \eqref{hp:cond1} is satisfied, $\lambda\in\Lambda$ and the conclusion follows directly from Theorem \ref{th:main1}.
\end{proof}

\begin{remark} \label{remark:f-nonneg-limiti}
	In Corollary \ref{corollary:lim1}, assumption \eqref{hp:limitea0} can be replaced by the weaker condition
	\begin{enumerate}[itemsep=0.2cm, label=$(\bar{\textnormal{h}}_\arabic*)$, ref = $\bar{\textnormal{h}_\arabic*}$]
		\setcounter{enumi}{6}
		\item\label{hp:limitea0-bar} $\limsup\limits_{t\to 0^+} \dfrac{F_{k_0}(t)}{t^2}=+\infty$
		\qquad for at least one $k_0 \in[1,N]$,
	\end{enumerate}
	without any sign condition on the functions $F_k$, by choosing in the proof a test function concentrated at $k_0$. If, in addition, the functions $f_k$ are nonnegative for all $k \in [1,N]$ and for all $t \ge 0$, then $\lambda^*_1=\lambda^*_2$, with
	\[
	\lambda^*_2 =  \dfrac{m^2}{2 a_2^2 s_+} \, \sup_{c>0}	\dfrac{c^2}{\displaystyle\sumo{N} F_k(c)}.
	\]
\end{remark}

\begin{corollary} \label{corollary:lim2}
	Let \eqref{H} be satisfied. Assume that
	\begin{enumerate}[itemsep=0.2cm, label=$(\widetilde{\textnormal{h}}_\arabic*)$, ref = $\widetilde{\textnormal{h}_\arabic*}$]
		\setcounter{enumi}{6}
		\item\label{hp:limitea0-f}	$\displaystyle\lim_{t\to0^+}\frac{f_k(t)}{t}=+\infty$ \quad for all  $k\in[1,N]$,
		\item\label{hp:limiteainfinito-f} $\displaystyle\lim_{t\to+\infty}\frac{f_k(t)}{t}=+\infty$  \quad for all  $k\in[1,N]$.
	\end{enumerate}
	Then, for every $\lambda \in (0, \lambda^*_1)$ problem \eqref{problem} admits at least two positive solutions  $u_1, u_2 \in X$ with opposite energy sign.
\end{corollary}
\begin{proof}
	Assumption \eqref{hp:limitea0-f} implies that $f_k$ is positive in a right neighborhood of $0$, so that $F_k\ge0$ for all $k\in[1,N]$ in a right neighborhood of $0$. Moreover, hypotheses \eqref{hp:limitea0-f}-\eqref{hp:limiteainfinito-f} imply the corresponding assumptions on the primitives \eqref{hp:limitea0}-\eqref{hp:limiteainfinito}. Therefore, the conclusion follows directly from Corollary \ref{corollary:lim1}.
\end{proof}

We conclude this section with two examples illustrating the applicability of our results. The first example is based on Corollary \ref{corollary:lim2} and involves a sign-changing potential term.

\begin{example}
	Let $N=2$, $s(k)=k+1$ for all $k\in[0,4]$, and
	\[
	q(k)=\frac{5k-6}{4}, \qquad k\in[1,2].
	\]
	Note that $q$ changes sign, $\mu_1=1$ and
	\[
	q_-=-\frac14>-\mu_1\,\frac{s_-}{s_+}=-\frac{s(1)}{s(3)}=-\frac12.
	\]
	For $t\ge0$, define
	\[
	f_1(t)=\sqrt t+t^3,
	\qquad
	f_2(t)=2\sqrt t+3t^3,
	\]
	and set $f_1(t)=f_2(t)=0$ for $t<0$. Hence, assumptions \eqref{hp:limitea0-f}-\eqref{hp:limiteainfinito-f} hold true. Moreover, the functions $f_1, f_2$ are nonnegative, so by Corollary \ref{corollary:lim2} and Remark \ref{remark:f-nonneg-limiti} it follows that the problem 
	\[
	\begin{cases}
		-\Delta^2u(k-1)- \dfrac{1}{k+2} \Delta u(k-1) + \dfrac{5k-6}{4} u(k)
		=\lambda f_k(u(k)), & k=1,2,\\
		u(0)=u(3)=0,
	\end{cases}
	\]
	admits at least two positive solutions with opposite energy sign for every $\lambda\in \left(	0,\, \frac1{12}\left(\frac12\right)^{1/5} \right)$.
\end{example}

Finally, the following example shows how Theorem \ref{th:main1} can be applied in the case of nonpositive coefficients.

\begin{example}
	Let $N=3$, $s(k)=k+1$ for all $k\in[0,5]$, and
	\[
	q(k)=-\frac{1}{10}, \qquad k\in[1,3].
	\]
	Then $q$ is nonpositive and
	\[
	q_-=-\frac{1}{10}>-\mu_1\,\frac{s_-}{s_+}=-\frac{2(2-\sqrt2)}{5}.
	\]
	Simple computations show that $\sigma=\frac{29}{5}$ and $K\approx0.079$.
	
	For $t\ge0$, define
	\begin{align*}
		f_1(t)&=\frac1{100}\bigl(t^{1/10}+t(1+\sin t)+t^3\bigr),\\
		f_2(t)&=\frac1{100}\bigl(t^{1/10}+t^2+t^3\bigr),\\
		f_3(t)&=\frac1{100}\bigl(t^{1/10}+e^t-1+t^3\bigr),
	\end{align*}
	and set $f_k(t)=f_k(0)=0$ for $t<0$. Then the functions $f_1,f_2,f_3$
	are continuous and nonnegative on $[0,+\infty)$. For $t\ge0$, their primitives are
	\begin{align*}
		F_1(t)& =\frac1{100}\left(\frac{10}{11}t^{11/10}
		+\frac12t^2+\sin t-t\cos t+\frac14t^4\right),\\
		F_2(t) &=\frac1{100}\left(\frac{10}{11}t^{11/10}
		+\frac13t^3+\frac14t^4\right),\\
		F_3(t) &=\frac1{100}\left(\frac{10}{11}t^{11/10}
		+e^t-t-1+\frac14t^4\right).
	\end{align*}
	Clearly, one has 
	\[
	L_\infty=+\infty.
	\]
	A direct computation shows that condition \eqref{hp:cond1} is satisfied by choosing $c=1$ and $d=\frac{1}{400}$, so that Theorem \ref{th:main1} applies. In this case, since $\Lambda=(0.242 , 2.151)$ approximately, we can choose $\lambda=1$ so that the problem
	\[
	\begin{cases}
		-\Delta^2u(k-1) - \dfrac{1}{k+2} \Delta u(k-1)-\dfrac{1}{10}u(k)
		= f_k(u(k)), & k\in[1,3],\\
		u(0)=u(4)=0,
	\end{cases}
	\]
	admits at least two positive solutions with opposite energy sign.
\end{example}

\section*{Acknowledgments}
\noindent The authors are members of the Gruppo Nazionale per l’Analisi Matematica, la Probabilità e le loro Applicazioni (GNAMPA) of the Istituto Nazionale di Alta Matematica (INdAM). This paper is supported by the project “An artificial intelligence approach for risk assessment and prevention of low back pain: towards precision spine care”, PNRR-MAD-2022-12376692, CUP: \mbox{J43C22001510001} funded by the European Union - Next Generation EU - NRRP M6C2 - Investment 2.1 Enhancement and strengthening of biomedical research in the NHS.


\end{document}